\documentclass[14pt]{mmcs_article}
\usepackage[english]{babel}
\usepackage{amsmath, amsthm, amsfonts, amssymb}

\usepackage[unicode,bookmarks,bookmarksopen,bookmarksopenlevel=0,colorlinks,%
pageanchor=true,breaklinks=true,linkcolor=blue,citecolor=blue,urlcolor=red,
pdfencoding = auto, psdextra]{hyperref}  

\newcommand{\rank}{\operatorname{rank}}

\usepackage{version}

\begin{document}

\addcontentsline{toc}{section}{Problem statement}

\section*{Calculation of chemical reactions in electrophoresis}

\begin{center}
\begin{doublespacing}
T. A. Petrukhina\footnote{Electronic address: petruhina@sfedu.ru} and
M. Yu. Zhukov\footnote{Electronic address: zhuk@math.sfedu.ru}

\textit{Institute of Mathematics, Mechanics and Computer Science,}

\textit{Southern Federal University, Rostov-on-Don, Russia}

(Dated: June, 2021)

\end{doublespacing}

\large{Abstract}
\end{center}

The main goal of the work is to find stationary solutions of the equations of chemical kinematics for mixtures of complex composition.

According to the given research goal the system of nonlinear equations for determining the equilibria of reversible chemical dissociation reactions was obtained and transformed in the work, an algorithm for solving the problem for tris-borate mixture was created and tested, allowing to calculate the pH of the solution taking into account the ionic force at given concentrations of the initial reagents of the mixture, as well as to process the results.

\textit{Keywords:} electrophoresis, isoelectrophocusing, tris-borate mixture.

\newpage
\addcontentsline{toc}{section}{Introduction}
\section*{Introduction}

The method of electrophoresis is one of the most effective methods for separating multicomponent chemically and biologically active mixtures into separate individual components by means of an electric field. The method is based on the ability of the mixture components (in particular, biopolymers or chemically active substances) to form charged complexes of molecules in solutions. The separation of the mixture into individual components occurs due to the fact that the migration rate in the electric field is individual for the individual mixture components. The movement under the action of electric field with different velocities eventually leads to the separation of the mixture.

There is a large number of literature on the practice and theory of electrophoresis, see e.g., \cite{Metody,Righetti1986,Righetti1986-1,Righetti1990,TroitskiyAzy,Shvedov,Bocek,ZhRStoy2001,BabskiiZhukovYudovichR,BabskZhukYud1989,
PolyakovaZhukovShiryaeva,Mosher1992,Zhukov2005,StepanovKorchemnaya,ZhukovShiryaeva2015MicroFluids,Osterman,Haglund1971,Everaerts1,Bier}. The mathematical theory of the electrophoresis method is developed in~\cite{BabskiiZhukovYudovichR,BabskZhukYud1989,PolyakovaZhukovShiryaeva,Zhukov2005}, see also~\cite{Mosher1992}. The method of electrophoresis is usually divided into three main types are the isoelectric focusing (IEF), zonal electrophoresis (ZE) and isotachophoresis (ITP). In fact, in practice, the three types of methods are subdivided into a large number of species, which, as a rule, differ from each other by different modifications, e.g., capillary zone electrophoresis, column isoelectrofocusing, etc. (for more information, see e.g., in~\cite{BabskiiZhukovYudovichR,BabskZhukYud1989,PolyakovaZhukovShiryaeva,Mosher1992}).

Let us focus on a slightly more detailed description of the isoelectric focusing method. It turns out that the rate of migration of the charged components of the mixture depends on the properties of the medium, more precisely on the  $\textrm{pH}$ medium (the concentration of hydrogen ions). Moreover, as a rule, biological substances (peptides, proteins, amino acids, etc.) are amphoteric compounds are the at some values of $\textrm{pH}$, they exhibit acidic properties, and at other values of basic properties. In other words, the individual components of the mixture can be both positively and negatively charged, depending on $\textrm{pH}$.
The migration rate in an electric field can be roughly represented in the form $v_k=q_k(pH) \gamma_k E$, where $v_k$ is the migration rate of the $k$ component, $q_k(pH)$ is the charge of the $k$ component, $\gamma_k$ is the characteristic mobility of the $k$ component (the rate per unit of the electric field strength; some constant depending, e.g., on the size of molecules), $E$ is the electric field strength fields. The function $q_k(pH)$ is an alternating sign. The value $pH=pI_k$, at which $q_k(pI_k) = 0$, is called the isoelectric point of the $k$ component. It is obvious that the component of the mixture at its isoelectric point is stationary.

\begin{figure}[H]
 \centering\includegraphics[scale=0.55]{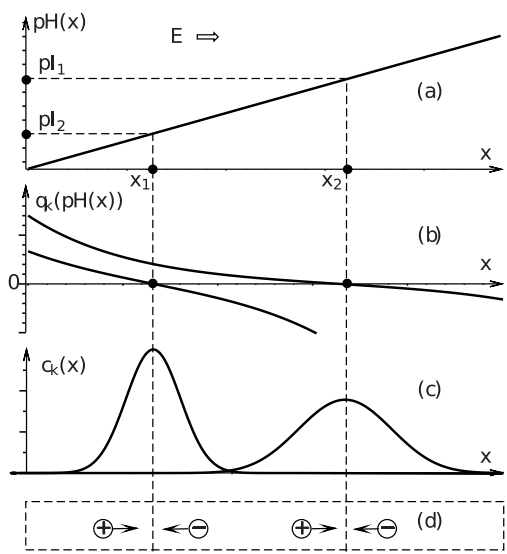}\\
 \caption{Isoelectric focusing scheme. Figure reproduced from~\cite[p.\,75]{PolyakovaZhukovShiryaeva}
with permission of the authors}
\label{ZS.fig:4.1}
\end{figure}

In isoelectric focusing, a part of the components of a multicomponent mixture ($S=(c_1,\dots,c_n)$) have small concentrations (substances to be separated), and the other part of the components ($B=(c_{n+1},\dots,c_N)$) has large concentrations (supporting or buffer mixture). Here $c_k$ are the components of the mixture (we also denote their concentrations at the same time), $N$ is the number of components of the mixture (excluding the solvent). The main role of the buffer mixture $B$ is to create an uneven distribution of $pH$ in the area in which the process takes place. In the spatially one-dimensional case, such us in a long cylindrical region or capillary, we can assume that $pH=pH(x)$, where $x$ is the coordinate. Another role of the buffer mixture $B$ to maintain the conductivity of the mixture. The substances of the buffer mixture $B$ are generally not amphoteric. In contrast, the substances released from the mixture, i.e., $S$ are amphoteric. If there are coordinates $x_k$ in the mixture in which $pH(x_k)=pI_k$, then it is obvious that, moving in an electric field, the components $S$ that hit the points $x_k$ will remain at these points. This is what the isoelectric focusing process is based on. Finally, at the given points $x_k$, the individual components of the mixture $c_k$ can be distinguished. The isoelectric focusing process is shown schematically in Figire~\ref{ZS.fig:4.1}.

Of course, there are many factors that should be taken into account when conducting (and mathematically modeling) the process. In particular, it is desirable to take into account the diffusion processes, the interaction between the components of the separated mixture $S$ and the components of the buffer mixture, etc.
It is also important to ensure as accurately as possible the distribution $pH$ of the solution over the region, i.e., setting $pH=pH(x)$. Note that, in principle, a situation is possible in which $pH(x) = \textrm{const}$. This case relates to zonal electrophoresis, in which the separation of the components of the mixture occurs not by isoelectric points, but by the rates of their migration.

\begin{figure}[H]
\centering\includegraphics[scale=0.75]{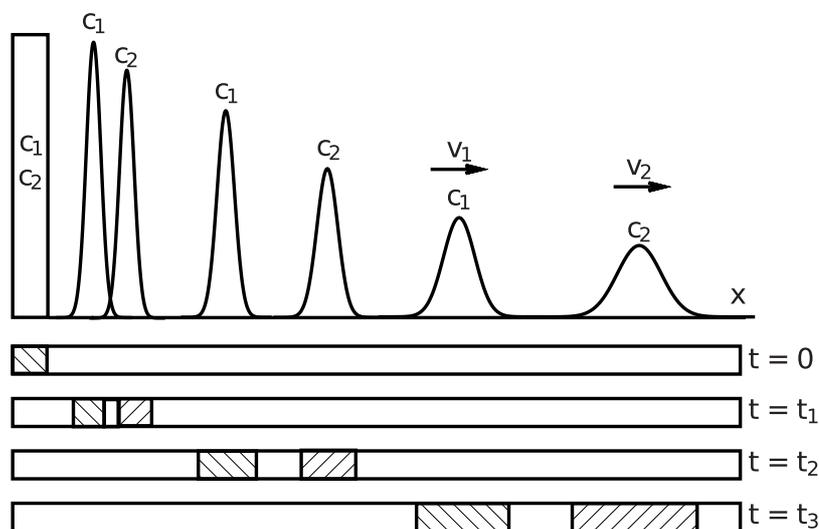} \caption{Scheme of separation of substances in zonal electrophoresis. Figure reproduced from~\cite[p.\,112]{PolyakovaZhukovShiryaeva} with permission of the authors}
\label{ZS.fig:5.1}
\end{figure}

An approximate diagram of the zonal electrophoresis process is shown in Figure~\ref{ZS.fig:5.1}.

As already mentioned, an important factor affecting the quality of the separation of the mixture is the accuracy of the $pH(x)$ function, which is mainly determined by the buffer mixture. As a rule, the concentrations of the components of the separated mixture $S$ are much smaller than the concentrations of the components of the buffer mixture $B$ and the influence of the components $S$ on $pH$ of the mixture as a whole is quite small. The value of $pH$ is formed as a result of chemical reactions of the components of the buffer mixture. Usually, the buffer mixture consists of a small amount of substances (two, three), which are acids and bases. As a result of dissociation reactions, acid and base ions are formed in the solution, as well as hydrogen ions $H^{+}$. The concentration of hydrogen ions is usually measured on a logarithmic scale. For these purposes, the ratio is used $pH = -\lg [H^{+}]$, where $H^{+}$ is the concentration of hydrogen ions, expressed in $\textrm{mol}/\textrm{L}$. As a rule, the rate of chemical reactions in the solution is quite high, and it can be assumed that the reactions occur ‘instantaneously’. In other words, the characteristic time of reaching the stationary state by chemical reactions is much less than the characteristic time of the transfer processes in electrophoresis. Thus, from a mathematical point of view, the determination of the $pH$ of a solution is reduced to finding the equilibria of the equations of chemical kinetics.

The main goal of the presented graduate work is precisely to find stationary solutions to the equations of chemical kinetics for mixtures of complex composition. The complex composition of a buffer mixture is a mixture in which dissociation reactions occur with the formation of a large number of different positive and negative ions. In contrast, an example of a mixture of a simple composition is, in particular, a mixture of one acid and one base, which dissociate into one positive ion and one negative ion. The problem of determining stationary solutions (equilibrium concentrations) is ultimately reduced to solving a rather complex system of nonlinear algebraic equations, the parameters of which are the dissociation constants of chemical reactions. In simple cases, when the dissociation constants are constant, it is possible in principle to reduce nonlinear algebraic equations to a system of linear equations, see e.g., \cite[p.\,351, 366--369]{Volpert}. However, even in this simplest case, the matrix of a system of linear equations is usually poorly conditioned, and obtaining an accurate numerical solution is a difficult task. At high concentrations of the mixture components, the dissociation ‘constants’ depend on the concentration values. In this case, the dissociation constants are usually called activity coefficients. This dependence makes it impossible to reduce a nonlinear system of equations to a linear one. In the case of dissociation reactions in aqueous solutions that occur with the formation of ions, the condition of the solution electroneutrality is added to the equations of chemical kinetics, which leads to an even greater complication of the problem of determining the equilibrium concentrations that satisfy the additional condition.

Due to the complexity of the problem, it is quite difficult to develop any universal algorithm suitable for any arbitrary reaction scheme. In this regard, the consideration is limited to the cases of buffer mixture of a fairly simple composition, which are usually used in electrophoresis. The dependence of the activity coefficients on the concentration of the components of the mixture is also limited to the case when the activity coefficients depend on some integral characteristic of the concentrations, and not on each concentration of the component separately. The ionic strength of the solution is chosen as such a characteristic. Some iterative algorithm together with also iterative algorithms for solving nonlinear systems of algebraic equations are proposed to solve the problem. As an example, a numerical solution is presented for some buffer mixture of a given composition. The effect of the ionic strength on the $pH$ of the solution is studied.

%
%
%
%

\newpage

\setcounter{equation}{0}
\setcounter{figure}{0}

\section{Basic equations}\label{Pt:sec:01}

We consider a buffer mixture consisting of $n$ components (excluding the solvent). We use the symbols $\xi_k$ ($k=0,\dots,n$) to denote the chemical components. The same designations are used for concentrations. As a rule, it is clear from the context what is being discussed, and cases of ambiguous interpretation are specified specifically. Usually, in chemistry, the symbol $A$ (or a set of symbols) is used to denote a substance, and the concentration of the substance is denoted by $[A]$. For example, the substance $NaCl$ has a concentration of $[NaCl]$. Further we use both forms of notation.

\subsection{Equations of chemical kinetics}\label{Pt:sec:01.01}

Consider a mixture in which $r$ chemical reactions of the form take place (see e.g., \cite[p.\,366--369]{Volpert}, \cite[p.\,52--56]{Zhukov2005})
\begin{equation}\label{Pt:01.01}
  \sum\limits_{k=0}^n \nu^+_{ik} \xi_k\
  \overset{k_i^+}{{\underset{k_i^-}{\rightleftharpoons}}}
  \sum\limits_{k=0}^n
  \nu^-_{ik} \xi_k, \quad i=1,\dots,r,
\end{equation}
where $\nu^+_{ik}$, $\nu^-_{ik}$ are the stoichiometric coefficients(non-negative integers); $k^{+}_i$, $k^{-}_i$ is the speed of direct and of the reverse reaction; $\xi_k$ are the name (symbols) chemicals.

The law of changes in concentrations $\xi_k(t)$ during chemical
reactions without taking into account mass transfer, i.e., the equations of chemical kinetics have the form
\begin{equation}\label{Pt:01.02}
  \frac{d \xi_k}{dt} =\sum\limits_{i=1}^{r} \nu_{ik}
     \sigma^{(i)} \equiv
     \sigma_k, \quad k=0,\dots,n,
\end{equation}
\begin{equation}\label{Pt:01.03}
  \sigma^{(i)} = -k^{+}_i \prod\limits_{k=0}^{n}
     \xi_k^{\nu_{ik}^{+}}
     +k^{-}_i \prod\limits_{k=0}^{n} \xi_k^{\nu_{ik}^{-}}, \quad
     \nu_{ik} = \nu_{ik}^{-} - \nu_{ik}^{+}.
\end{equation}
Here $\xi_k$ are the molar concentrations of the components ($\textrm{mol}/\textrm{L}$),
$\sigma_k$ is the density of mass sources (concentrations). In chemistry, the values of $\sigma_k$ are unfortunate to be called chemical reaction rates, although, by implication, they are more likely to be concentration change rates.

\subsection{Equilibrium conditions in the general case}\label{Pt:sec:01.02}

The necessary and sufficient condition for the existence of an equilibrium, i.e., stationary solution of the system \eqref{Pt:01.02} is the solvability
of the system \quad $\sigma^{(i)} = 0$, $i=1,\dots,r$, which after
logarithm is reduced to a system of equations with respect to $\ln \xi_k$
\begin{equation}\label{Pt:01.04}
   \sum\limits_{k=0}^{n} \nu_{ik} \ln \xi_k =\ln K_i, \quad
     K_i=\frac{k_i^+}{k_i^-}, \quad i=1,\dots,r.
\end{equation}
Here $K_i$ are the equilibrium constants of chemical reactions (the dissociation constants in the case of reactions with the formation of ions).

In the simplest cases, the values of $K_i$ are constant (and the name \emph{constant} is justified). In general, the values of $K_i$ depend on the concentrations, but the name of the ‘constant’ is traditionally preserved.

A more correct notation of the equations \eqref{Pt:01.04}, which emphasizes the dependence of $K_i$ on concentrations, has the form
\begin{equation}\label{Pt:01.05}
   \sum\limits_{k=0}^{n} \nu_{ik} \ln \xi_k =\ln K_i(\xi), \quad
     \xi =(\xi_1, \dots, \xi_n).
\end{equation}

\subsection{Linear equations for $K_i=\textrm{const}$}\label{Pt:sec:01.03}

In the case where the constants $K_i$ are constant, the system \eqref{Pt:01.04}
\begin{equation}\label{Pt:01.06}
   \sum\limits_{k=0}^{n} \nu_{ik} \ln \xi_k =\ln K_i, \quad
    i=1,\dots,r.
\end{equation}
the system is a \emph{linear} system of equations with respect to the unknowns $\ln \xi_k$.

Let the rank of the matrix $(\nu_{ik})_{i=1,\dots,r}^{k=0,\dots,n}$ be equal to $r_0=n-r$. Then there is a fundamental system of solutions
$\lambda_k^{(s)}$, $s=1,\dots,r_0$ of a homogeneous
equation~\eqref{Pt:01.06}
\begin{equation}\label{Pt:01.07}
  \sum\limits_{k=0}^{n} \nu_{ik} \lambda_k^{(s)} =0,\quad
      s=1,\dots,r_0,\quad  i=1,\dots,r.
\end{equation}

Obviously, in this case there are independent linear first
integrals of the equations \eqref{Pt:01.02}
\begin{equation}\label{Pt:01.08}
 a_s=\sum\limits_{k=0}^{n} \lambda_k^{(s)} \xi_k, \quad
     \frac{da_s}{dt}=0, \quad  s=1,\dots,r_0.
\end{equation}
In chemistry, the values of $a_s$ are called analytical concentrations.

In~\cite[p.\,52--56]{Zhukov2005} the concept of a $s$ is the chemical subsystem is introduced (see also
\cite{BabskiiZhukovYudovichR}, p.~45, where a similar
definition). The set of concentrations $\xi_k$ forms the $s$ is the chemical subsystem $A_s$ if $\lambda_k^{(s)}
\ne 0$, i.e.,
\begin{equation}\label{Pt:01.09}
   A_s=
     \left\{ \xi_k \colon \lambda_k^{(s)} \ne 0, \quad k=0,\dots,n.
    \right\}.
\end{equation}
The values of $a_s$  are called concentrations of the chemical subsystem. Note that the same concentrations of $\xi_k$ can belong to different chemical subsystems.

\subsection{Mixture with charged components}\label{Pt:sec:01.04}

Since the analytical concentrations of $a_s$ are determined based on the homogeneous equations \eqref{Pt:01.06}, the same procedure can be used for the nonlinear equations \eqref{Pt:01.05}. In other words, the allocation of analytical concentrations for the systems  \eqref{Pt:01.05}, \eqref{Pt:01.06} is the same.
Next, we consider a more general system \eqref{Pt:01.05}.

In the case of a buffer mixture with charged components, the condition of electroneutrality must be added to the equations \eqref{Pt:01.06}
\begin{equation}\label{Pt:01.10}
   \sum\limits_{k=0}^{n}z_k\xi_k = 0,
\end{equation}
where $z_k$ are the charges, more precisely charges, components (charge in units of electron charge, integers).

In essence, it is another chemical subsystem is the a set of charged components
\begin{equation}\label{Pt:01.11}
   A= \left\{ \xi_k{:} \ \lambda_k=z_k \ne 0, \quad
  k=0,\dots,n \right\}.
\end{equation}
The analytical concentration of such a subsystem is the molar charge of the solution equal to zero.

The presence of the electroneutrality equation makes it possible to express one of the concentrations  through all the others. Traditionally, the concentration of hydrogen ions or the acidity of the solution is chosen as such a variable in aqueous solutions
\begin{equation}\label{Pt:01.12}
   pH=- \lg H^{+}, \quad \xi_0 = H^{+},
\end{equation}
where the concentration of $H^{+}$ is expressed in $\textrm{mol}/\textrm{L}$.

Adding the electroneutrality equation \eqref{Pt:01.10} to \eqref{Pt:01.05} significantly complicates the general system of equations. Indeed,
even in the case when $K_i$ is constant, the equations \eqref{Pt:01.05} are linear on the variables $\ln\xi_i$, while the equation \eqref{Pt:01.10} is linear on the variables $\xi_i$.

Some simplification of the problem can be achieved by introducing for each chemical subsystem the degrees of dissociation is the relative fraction of the concentration of subsystem components from the analytical concentration of the subsystem
\begin{equation}\label{Pt:01.13}
  \xi_k=\theta_k^{s}(\xi_n, a_1,\dots,a_{r_0-1}) a_s, \quad
     k=1,\dots,n-1.
\end{equation}
Here $\theta_k^{s}$ are the degrees of dissociation of the $k$ component for $s$ the chemical subsystem, i.e., the concentration fraction of the $k$ component in the concentration of $a_s$.

Obviously, (see \eqref{Pt:01.09})
\begin{equation}\label{Pt:01.14}
   \sum\limits_{k=1}^{n} \theta_k^s \lambda_k^{(s)} =0, \quad
     s=1,\dots,r_0, \quad r_0=n-r.
\end{equation}

The introduction of dissociation degrees allows to write the system \eqref{Pt:01.05}, \eqref{Pt:01.10} only in terms of analytical $a_a$ and concentration $H^{+}$, which greatly simplifies the system, reducing the number of unknowns. Note that in the simplest cases of electrophoresis, the degrees of dissociation $\theta_k^s$ depend only on the concentration of hydrogen ions $H^{+}$ or $pH$ of the medium.

\subsection{Example. Mixture of acid and base}\label{Pt:sec:01.05}

Consider a buffer mixture consisting of a single-charge acid and a single-charge base. In the description we use chemical notation, which in practical cases is more convenient than formalized (faceless) mathematical notation.

The scheme of acid and base dissociation reactions has the form (compare with \eqref{Pt:01.01}, and also see e.g., \cite[p.\, 61, 62]{PolyakovaZhukovShiryaeva})
\begin{equation}\label{Pt:01.15}
  HA   \overset{K_a}\rightleftarrows A^- + H^+,\quad
  H^+B \overset{K_b}\rightleftarrows B + H^+,
\end{equation}
where $HA$ it the acid, $A^-$ is the acid residue, $B$ is the base, $H^+B$ is the base residue,
$K_a$, $K_b$ is the dissociation constants of the acid and base respectively.

Introduce the molar analytical concentrations of acid $a$, base
$b$ and degrees of dissociation is the $\alpha$ for the acid ion and $\beta$ for the base ion (compare~with~\eqref{Pt:01.08})
\begin{equation}\label{Pt:01.16}
  a=[HA]+[A^-],\quad
  b=[H^+B]+[B].
\end{equation}
A well-known mnemonic rule is used here is the sum of the concentrations of substances that have the same symbol in their names is the analytical concentration of the substance. In this case, the analytical concentration of the acid $a$ consists of the concentration of the acid ion $A^-$ and the neutral acid residue $HA$.
The analytical concentration of the base $b$ consists of the concentration of the base ion $H^+B$ and the neutral residue of the base $B$. Note that the analytical concentration is a real measurable value in contrast to the concentration of ions and the concentration of the residue. In other words, we can say that the solution contains so much acid ($a$), but it is impossible (without using special methods or calculations) to say how many ions or neutral molecules formed during dissociation are contained.

Let us introduce the degrees of dissociation (compare~with~\eqref{Pt:01.13})
\begin{equation}\label{Pt:01.17}
  \alpha = \dfrac{[A^-]}{a},\quad \beta = \dfrac{[H^+B]}{b}.
\end{equation}

The equilibrium conditions of the reactions \eqref{Pt:01.15} have the form (compare with \eqref{Pt:01.03}--\eqref{Pt:01.06})
\begin{equation}\label{Pt:01.18}
  \dfrac{[A^-][H^+]}{[HA]}=K_a, \quad
  \dfrac{[B][H^+]}{[H^+B]}=K_b.
\end{equation}
Then, taking into account  \eqref{Pt:01.17}, we obtain expressions for the degrees of dissociation, which depend only on the concentration of hydrogen ions (in the case under consideration)
\begin{equation}\label{Pt:01.19}
  \alpha=\dfrac{K_a}{K_a+[H^+]},\quad \beta=\dfrac{[H^+]}{K_b+[H^+]}.
\end{equation}

The electroneutrality equation is used to determine the value of $[H^+]$ (compare with \eqref{Pt:01.10})
\begin{equation}\label{Pt:01.20}
[H^+B] - [A^-] + [H^+] - [OH^-] = 0.
\end{equation}
Here $[OH^-]$ is the concentration of hydroxyl ions.

Recall that the hydroxyl ions  $OH^-$ arise as a result of the reaction of autoprotolysis of water
\begin{equation}\label{Pt:01.21}
  H_2O \overset{K^2_w}\rightleftarrows OH^- + H^+,\quad
\end{equation}
where $H_2O$ is the water, $OH^-$ is the hydroxyl ion, $H^+$ is the hydrogen ion, $K^2_w$ are the water autoprotolysis (dissociation) constants.

The reaction equilibrium condition \eqref{Pt:01.21} has the form
\begin{equation}\label{Pt:01.22}
  [OH^-][H^+] = K^2_w \quad \Rightarrow \quad [OH^-] = \dfrac{K^2_w}{[H^+]}.
\end{equation}
In practice, the concentration of $[OH^-]$ is much less than the sum (algebraic) of the remaining terms of the equation \eqref{Pt:01.20} and in calculations it is often neglected, writing the equation of electroneutrality \eqref{Pt:01.20} in a simplified form
\begin{equation}\label{Pt:01.23}
[H^+B] - [A^-] + [H^+] = 0.
\end{equation}
In fact, the concentration of $[H^+]$ is also small, and one can further simplify the equation by writing
\begin{equation}\label{Pt:01.24}
[H^+B] - [A^-] = 0.
\end{equation}

Using \eqref{Pt:01.17}, \eqref{Pt:01.19} and \eqref{Pt:01.24}, getting
\begin{equation*}
\beta b - \alpha a = 0
\end{equation*}
or
\begin{equation}\label{Pt:01.25}
  \dfrac{K_a}{K_a+[H^+]}a=\dfrac{[H^+]}{K_b+[H^+]}b.
\end{equation}
This equation is a square equation with respect to $[H^+ ]$, the positive solution of which determines
the value of the $pH$ of the buffer mixture
\begin{equation}\label{Pt:01.26}
  [H^+]=\dfrac12 K_a
  \left(\dfrac{a}{b}-1+\sqrt{\left(\dfrac{a}{b}-1\right)^2+
  4\dfrac{K_b}{K_a}\dfrac{a}{b}}
  \right),
\end{equation}
\begin{equation*}
  pH=-\lg[H^+].
\end{equation*}
Various approximate expressions for the ratio \eqref{Pt:01.26}, made on the basis of assumptions about the orders of magnitude of concentrations and dissociation constants are known as the Henderson—Hasselbalch relations. As a matter of fact, all simplifications of the equation of electroneutrality
\eqref{Pt:01.20} before \eqref{Pt:01.23} and then \eqref{Pt:01.24} are associated with the possibility of obtaining a simple solution of the form \eqref{Pt:01.26}.
It is interesting to note that the value of $pH$  depends only on the ratio of the concentrations of $a/b$ and the dissociation constants of $K_b/K_a$.

\begin{figure}[H]
  \centering
  \includegraphics[scale=0.75]{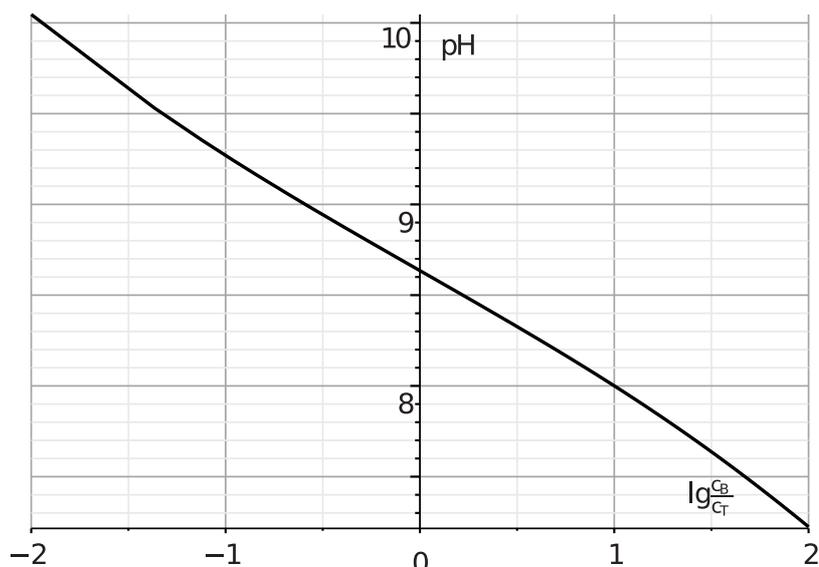}
  \caption[Dependence of $pH$ on $\lg(C_B/C_T)$ for the Tris-borate buffer]{Dependence of $pH$ on $\lg(c_B/c_T)$  for the Tris-borate buffer, calculated by the formula \eqref{Pt:01.26}, in the case of $\textrm{pK}_B = 9.29$, $\textrm{pK}_T = 7.98$. Here $C_B=a$, $C_T=b$, $Ka = -\lg\textrm{pK}_B$, $Kb = -\lg\textrm{pK}_T$. Figure is reproduced from \cite[p.\,62]{PolyakovaZhukovShiryaeva} with permission of the authors}
  \label{ZS.fig:3.1}
\end{figure}

The results of the calculation using the formula \eqref{Pt:01.26} do not always coincide with the results of the experiment. There are quite a lot of reasons for this. One of them is the dependence of the dissociation constant on the concentration, temperature, the  performed simplifications of the electroneutrality equation, etc. In addition, different literature sources provide different data for the dissociation constants. For example, in Figure~\ref{ZS.fig:3.1} shows the calculation results for a Tris-borate buffer consisting of boric acid and tris (base) are given. Recall the chemical formulas: boric acid $H_3BO_3$, tris $C_4H_{11}NO_3$ (2-amino-2-hydroxymethyl-propane-1, 3-diol).

\subsection{Influence of the ionic strength of the solution}\label{Pt:sec:01.06}

The formula \eqref{Pt:01.26} is obtained under the assumption that the values $K_a$, $K_b$ are constant.
In reality, the values $K_a$, $K_b$ depend on the concentrations of the components of the mixture. There are a large number of theories and semi-empirical formulas that allow us to take into account this dependence. One of the simplest variants is the hypothesis that the dependence of the dissociation constants is determined by some integral characteristic of the mixture, namely, the ionic strength of the solution $I$ (see e.g., \cite[p.\,62--64]{PolyakovaZhukovShiryaeva}, and also \cite[p.\;41]{Korita})
\begin{equation}\label{Pt:01.27}
  I = \frac{1}{2}\sum\limits_{k}\xi_i z_i^2,
\end{equation}
where $\xi_i$ is the molar concentration of ions, $z_i$ are the charges of ions.

Note that $I$ differ in different schools according to the scales in which the ion concentrations are expressed. If, as usual, the ionic strength is expressed in $\textrm{mol}/\textrm{L}$, then (see e.g., \cite[p.\; 42, 43, 51]{Korita})
\begin{equation}\label{Pt:01.28}
  \lg \gamma_k = - A z_k^2 I^{1/2}, \quad A = 0.5093 ~\text{mol}^{-1/2}\text{L}^{1/2},
\end{equation}
where $\gamma_k$ are the activity coefficients. These ratios are obtained for strong electrolytes and there are many refinements of the formula, usually for specific reactions. Taking into account the activity coefficients, the equilibrium conditions \eqref{Pt:01.18} should be written in the form
\begin{equation}\label{Pt:01.29}
  \dfrac{[A^-][H^+]}{[HA]}\dfrac{\gamma_{A^-} \gamma_{H^+}}{\gamma_{HA}} = K_a, \quad
  \dfrac{[B][H^+]}{[H^+B]} \dfrac{\gamma_B \gamma_{H^+}}{\gamma_{H^+B}}= K_b.
\end{equation}
In other words, in the equilibrium conditions \eqref{Pt:01.05} one should make substitutions
\begin{equation}\label{Pt:01.30}
  \xi_k \to b_k = \gamma_k \xi_k.
\end{equation}
The values of $b_k$ are called concentration activities.
The activity coefficients of neutral components are considered equal to one ($\gamma_k =1$, if $z_k = 0$). In the case of \eqref{Pt:01.29}, the activity coefficients of the neutral components are considered equal to one, i.e., $\gamma_{HA} = 1$, $\gamma_B = 1$.

When speaking about the influence of the ionic strength on the dissociation constants, we mean rewriting the relations \eqref{Pt:01.29} in the form
\begin{equation}\label{Pt:01.31}
  \dfrac{[A^-][H^+]}{[HA]} = K'_a = \dfrac{\gamma_{HA}}{\gamma_{A^-} \gamma_{H^+}} K_a, \quad
  \dfrac{[B][H^+]}{[H^+B]} = K'_b = \dfrac{\gamma_{H^+B}}{\gamma_B \gamma_{H^+}}K_b.
\end{equation}
In other words, new dissociation constants  $K'_a$, $K'_b$, are introduced, which, formally, allows us to preserve the previous form of the equations \eqref{Pt:01.05}.
Of course, this does not lead to any simplifications, since the new ‘constants’ depend on the ionic strength and, therefore, by virtue of \eqref{Pt:01.27}, \eqref{Pt:01.28} on the concentrations of the solution components.

The effect of the ionic strength on the calculation of the $pH$ of the solution for the reactions \eqref{Pt:01.15}  is shown in Table \ref{ZS.Tab:8}.
\begin{table}[H]
\caption[Effect of ionic strength on $pH$ and dissociation constants]{Effect of ionic strength on $pH$ and dissociation constants $\textrm{pK}_B$, $\textrm{pK}_T$. The lower three lines correspond to $I = 0$. The activity coefficients of single-charge ions coincide with each other and are indicated in the table as~$\gamma$. The table is reproduced from \cite[p.\,63]{PolyakovaZhukovShiryaeva} with permission of the authors}
\label{ZS.Tab:8}

 \centering
\begin{tabular}{|l|l|l|l|}
  \hline
                   & $a = 0.1$ & $a = 0.1$ & $a = 0.2$ \\
                   & $b = 0.2$ & $b = 0.1$ & $b = 0.1$ \\
  \hline
  $I$              & $0.029653              $ & $0.020767                 $ & $0.029650                 $ \\
  $\gamma$         & $0.817142              $ & $0.844513                 $ & $0.817151                 $ \\
  $\textrm{pK}'_a$  & $9.114595              $ & $9.143213                 $ & $9.114605                $ \\
  $\textrm{pK}'_b$  & $7.980000              $ & $7.980000                 $ & $7.980000                $ \\
  $\textrm{pH}$    & $\textbf{8.739}        $ & $\textbf{8.562}           $ & $\textbf{8.355}           $ \\
  \hline
  $\textrm{pK}_a$  & $9.29                  $ & $9.29                     $ & $9.29                     $ \\
  $\textrm{pK}_b$  & $7.98                  $ & $7.98                     $ & $7.98                     $ \\
  $\textrm{pH}$    & $\textbf{8.819}        $ & $\textbf{8.635}           $ & $\textbf{8.451}           $ \\

  \hline
\end{tabular}
\end{table}
It should be said that saving a large number of characters in the table is done in order to see the differences in the results. In practice, it is enough to be limited to two or three decimal places. The purpose of the entire section~\ref{Pt:sec:01} was to show different approaches (abstract and concrete) to solving the problem of determining the $pH$ of a buffer mixture and to demonstrate the complexities of the equations and ways to simplify them.

\newpage
\setcounter{equation}{0}
\setcounter{figure}{0}

\section{Problem statement and solution algorithm}\label{Pt:sec:02}

The main purpose of the presented work is to calculate the chemical equilibrium and determine the $pH$ of the buffer mixture for the Tris-borate buffer, which is often used in the processes of isoelectrofocusing and zonal electrophoresis. Despite the names that match the example of par.~\ref{Pt:sec:01.05} and the fact that the mixture consists of boric acid and tris (base), the actual dissociation scheme is much more complex than \eqref{Pt:01.15}. The fact is that boron compounds
are able to create complex charged complexes in solutions, and dissociation occurs in a solution in a rather complex way. In addition, as the available experimental data show, the ionic strength of the solution has a significant effect on the $pH$  of the mixture (see~\eqref{Pt:01.15}).

\subsection{Problem statement}\label{Pt:sec:02.01}

For a system of chemical reactions \eqref{Pt:01.01} describing dissociation processes (the solution contains charged components)
\begin{equation}\label{Pt:02.01}
  \sum\limits_{k=0}^n \nu^+_{ik} \xi_k\
  \overset{k_i^+}{{\underset{k_i^-}{\rightleftharpoons}}}
  \sum\limits_{k=0}^n
  \nu^-_{ik} \xi_k, \quad i=1,\dots,r,
\end{equation}
is necessary to determine the equilibrium concentrations, taking into account the electroneutrality of the mixture, in the case when the dissociation constants depend only on the ionic strength of the solution.

From a mathematical point of view, the problem is to find a solution to a system of nonlinear equations
(see~\eqref{Pt:01.05}, \eqref{Pt:01.10}, \eqref{Pt:01.27}, \eqref{Pt:01.28})
\begin{equation}\label{Pt:02.02}
   \sum\limits_{k=0}^{n} \nu_{ik} \ln \xi_k = \ln K'_i(I), \quad
     i = 1,\dots,r,
\end{equation}
\begin{equation}\label{Pt:02.03}
   \sum\limits_{k=0}^{n}z_k\xi_k = 0,
\end{equation}
where
\begin{equation}\label{Pt:02.04}
    \ln K'_i(I) = \ln K_i - \sum\limits_{k=0}^{n} \nu_{ik} \ln \gamma_k,
\end{equation}
\begin{equation}\label{Pt:02.05}
  \ln \gamma_k = - A z_k^2 I^{1/2}\ln10,
\end{equation}
\begin{equation}\label{Pt:02.06}
  I = \frac{1}{2}\sum\limits_{k=0}^{n}\xi_i z_i^2,
\end{equation}
\begin{equation}\label{Pt:02.07}
  \rank \nu_{ik} = r_0 = n -r > 0.
\end{equation}
Parameters (constant) are considered to be the specified values $\nu_{ik}$, $K_i$, $z_k$, $A$.

Equations \eqref{Pt:02.02} are the equilibrium equations of chemical reactions, the equation \eqref{Pt:02.03} is the equation of electroneutrality. The relations \eqref{Pt:02.04}--\eqref{Pt:02.06} are defining relations that define the dependence of $K'_i$ on $I$, which in turn depends on $\xi_k$.

The unknown quantities to be determined are the concentrations $\xi_k$. By virtue of the \eqref{Pt:02.07} the system \eqref{Pt:02.02} (and the entire system \eqref{Pt:02.02}--\eqref{Pt:02.06} in general) is undefined. The number of unknowns $\xi_k$ is $(n+1)$, while the number of equations is $(r+1)$.

The system of equations \eqref{Pt:02.02} for constants $K'_i$ is a system of linear equations with respect to the variables $\ln\xi_k$. Due to the condition \eqref{Pt:02.07} this allows us to determine the values of $a_s$ (analytical concentrations) such that the relations are satisfied (see \eqref{Pt:01.07}, \eqref{Pt:01.08})
\begin{equation}\label{Pt:02.08}
 a_s=\sum\limits_{k=0}^{n} \lambda_k^{(s)} \xi_k, \quad  s=1,\dots,r_0.
\end{equation}
where $\lambda_k^{(s)}$, $s=1,\dots,r_0$ is the fundamental solution of a system of homogeneous equations
\begin{equation}\label{Pt:02.09}
  \sum\limits_{k=0}^{n} \nu_{ik} \lambda_k^{(s)} =0,\quad
      s=1,\dots,r_0 = n-r, \quad  i=1,\dots,r.
\end{equation}
Thus, if to equations \eqref{Pt:02.02}, \eqref{Pt:02.03}, considering $a_s$ known, add equations \eqref{Pt:02.08}, the number of equations to determine $\xi_k$ will match the number of unknowns.

\textbf{Final statement of the problem} is the at given parameters $A$, $\nu_{ik}$, $K_i$, $z_k$,  $a_s$
($i=1,\dots,r$, $k=0,\dots,n$, $s=1,\dots,r_0=n-r$) it is required to find a solution to the equations \eqref{Pt:02.02}--\eqref{Pt:02.09} regarding unknowns $\xi_k$, ($k=0,\dots,n$).

\textbf{Remark 1.} The equation \eqref{Pt:02.02} can be replaced by the equation \label{Pt:note:01}
\begin{equation}\label{Pt:02.10}
   \prod\limits_{k=0}^{n} \xi_k^{\nu_{ik}} = K'_i(I).
\end{equation}
This is convenient to do when constructing analytical transformations. On the contrary, when using numerical algorithms.

\subsection{Solution algorithm}\label{Pt:sec:02.02}

In principle, the system of nonlinear equations \eqref{Pt:02.02}--\eqref{Pt:02.09} can be solved by ordinary numerical methods, for example, the Newton's method (see e.g., \cite{Kalitkin,Bahvalov}).
However, as is well known, such a method requires good initial approximations and may be inefficient for large system dimensions (for large $n$). One of the reasons for the inefficiency is the large ‘spread’ of the values of the $K_i$ parameters. In the problems of electrophoresis, as a rule, $K_i = 10^{-pK_i}$, where
$2 < pK_i < 12$. Moreover, in the physical sense, concentrations must obviously satisfy the inequalities $\xi_k \geq 0$. In fact, the inequalities must be strict $\xi_k > 0$. Cases where one (or more) concentrations turn are degenerate and should be considered separately.

The problem with finding only positive solutions is solved quite easily. It is sufficient to substitute variables in the equations
\begin{equation}\label{Pt:02.11}
   \xi_k = e^{u_k} > 0, \quad u_k = \ln\xi_k.
\end{equation}
Note that with this substitution, the equations \eqref{Pt:02.02} for constants $K'_i$ are linear equations with respect to the new unknowns $u_i$. In contrast, the equations \eqref{Pt:02.03}, \eqref{Pt:02.08}, and the relation \eqref{Pt:02.06}, being linear in the variables $\xi_k$, will include nonlinear terms of the form $e^{u_k}$.  An additional advantage of replacing \eqref{Pt:02.11} is the ‘normalization’ of the parameters $\ln K'_i$, which for the problems of electrophoresis will be of the order of several units.

In principle, it is possible to create some universal algorithm for solving nonlinear equations \eqref{Pt:02.02}--\eqref{Pt:02.09}, solving the problem with the choice of the initial approximation, for example, by means of a parameter motion. However, this use of the algorithm will not be very convenient. Firstly, the equations
\eqref{Pt:02.02}--\eqref{Pt:02.09} will contain a large number of parameters, many of which are zero. This applies to the matrix of stoichiometric coefficients $\nu_{ik}$ and the charges $z_k$. In other words, a large number of calculations will be done ‘for nothing’. Secondly, the definition of the fundamental solutions of the equations \eqref{Pt:02.09} is obviously not unambiguous. In particular, such solutions are defined with multiplier accuracy, which may lead to unreasonable expressions for analytical concentrations of $a_s$. In particular, concentrations that are not always clearly interpreted from a chemical point of view can be obtained.

Note that the creation of a general algorithm is not the goal of the graduate work, which is focused on solving the problem of specific chemical reactions.

Let us formulate several steps  to solve the problem \eqref{Pt:02.02}--\eqref{Pt:02.09}, which significantly facilitate the construction of the solution.

\label{Pt:etap:01}
\textbf{Step 1.} Find the fundamental solutions $\lambda_k^{(s)}$, $s=1,\dots,r_0$ of the homogeneous system of equations (see \eqref{Pt:01.07} and \eqref{Pt:01.09})
\begin{equation}\label{Pt:02.12}
  \sum\limits_{k=0}^{n} \nu_{ik} \lambda_k^{(s)} =0,\quad
      s=1,\dots,r_0,\quad  i=1,\dots,r.
\end{equation}
This is done by using  $\lambda_k^{(s)}$  isolate the chemical subsystems of  $A_s$ and obtain ratios for the concentrations of $a_s$ (see \eqref{Pt:01.07}, \eqref{Pt:01.08}, \eqref{Pt:02.08})
\begin{equation}\label{Pt:02.13}
 a_s=\sum\limits_{k=0}^{n} \lambda_k^{(s)} \xi_k, \quad
 s=1,\dots,r_0.
\end{equation}
\begin{equation}\label{Pt:02.14}
   A_s=
     \left\{ \xi_k \colon \lambda_k^{(s)} \ne 0, \quad k=0,\dots,n.
    \right\},
\end{equation}
In doing so, one should strive for ‘chemical clarity’\ of the result. It is desirable that $a_s$ have a clear chemical meaning.

\label{Pt:etap:02}
\textbf{Step 2.} Assuming $K'_i$ to be constant, we solve the equations  \eqref{Pt:02.02}, replacing them with  \eqref{Pt:02.10}, together with the equations  \eqref{Pt:02.08} (or
\eqref{Pt:02.13}). Here one should try to construct a solution of the form
\begin{equation}\label{Pt:02.15}
 \xi_k=\xi_k(\xi_0,a_1,\dots, a_{r_0}; K'_1,\dots,K'_n), \quad k = 1,\dots,n. 
\end{equation}
In other words, we should express the concentrations of all components  through the analytical concentrations $(a_1,\dots, a_{r_0})$, the dissociation constants $(K'_1,\dots,K'_n)$ and the concentration of $\xi_0$. Here the index $k=0$ corresponds to the ion $H^+$.

Note that, given \eqref{Pt:02.04},   \eqref{Pt:02.05} the dependency \eqref{Pt:02.15} is representable in form
\begin{equation}\label{Pt:02.16}
 \xi_k=\xi_k(\xi_0,a_1,\dots, a_{r_0}, I), \quad k =1,\dots,n. 
\end{equation}

If it is impossible to construct an explicit solution of the form \eqref{Pt:02.16} we create a numerical algorithm that allows us to construct functions $\xi_k$ for all components.
The input parameters of such an algorithm should be the values $(\xi_0,a_1,\dots, a_{r_0}, I)$, and the output parameters of the concentration $\xi_k$ (except, of course, $\xi_0$).

Next, consider that the concentrations of $\xi_k$ are determined in one way or another (i.e., explicitly or numerically) by the relations \eqref{Pt:02.16}.

\label{Pt:etap:03}
\textbf{Step 3.} Using \eqref{Pt:02.16}, we substitute the concentrations of the charged components in the equation of electroneutrality \eqref{Pt:02.03}
\begin{equation}\label{Pt:02.17}
   q(\xi_0,I,a_1,\dots, a_{r_0})\equiv z_0\xi_0 + \sum\limits_{k=1, z_k \ne 0}^{n}z_k\xi_k(\xi_0,a_1,\dots, a_{r_0}, I) =0.
\end{equation}
In addition, given \eqref{Pt:02.06}, we have
\begin{equation}\label{Pt:02.18}
  \varphi(\xi_0,I,a_1,\dots, a_{r_0})\equiv I - \frac{1}{2}z_0^2\xi_0 - \frac{1}{2}\sum\limits_{k=1, z_k \ne 0}^{n} z_i^2\xi_k(\xi_0,a_1,\dots, a_{r_0}, I)=0.
\end{equation}
Thus, a system of two nonlinear equations of the form (parameters $(a_1,\dots, a_{r_0})$ omitted)
\begin{equation}\label{Pt:02.19}
   q(\xi_0,I) = 0, \quad \varphi(\xi_0,I) = 0,
\end{equation}
the solution of which determines the concentration of hydrogen ions $\xi_0=H^+$ (and thus the $pH$ of the buffer mixture) and the ionic strength of the solution $I$. The initial approximation, in particular for the Newton's method, can be obtained by constructing, for example, level lines for the functions $q(\xi_0,I)$, $\varphi(\xi_0,I)$.

In the case where the influence of the ionic strength does not affect the dissociation constants, in the relations \eqref{Pt:02.16} we should put $I=0$ and consider only
one equation of the system \eqref{Pt:02.19}
\begin{equation}\label{Pt:02.20}
   q(\xi_0,I)\bigr|_{I=0} = 0.
\end{equation}

\textbf{Remark 2.} If information about the concentration of all the components of the solution is not required, then it is sufficient to limit the construction of the dependences \eqref{Pt:02.16} only for charged components, since only such components participate in the equation of electroneutrality \eqref{Pt:02.03} and in the ratio \eqref{Pt:02.06} for the ionic strength. Moreover, it may turn out that the system \eqref{Pt:02.10}, \eqref{Pt:02.13} is easier to solve not for the set $(\xi_, a_1,\dots, a_{r_0})$, and for any other set of ‘parameters’. It is important that such a set consists of $\xi_0$ and any $r_0$ values, with respect to which the system is easily solved. The main goal of Step $2$ is to reduce the dimensionality of the system of initial nonlinear equations.  If the transformation is successful, the dimensionality of the original system ($n+1$ with constants $K'_i$) is reduced to $r_0+1 = n-r +1$.
\label{Pt:note:02}

\newpage
\setcounter{equation}{0}
\setcounter{figure}{0}

\section{Calculation of $pH$ for Tris-borate buffer mixture}\label{Pt:sec:03}

Let us give an example of using the results described in par.~\ref{Pt:sec:02.02}. For these purposes, consider the complete dissociation scheme in an aqueous solution of boric acid + tris. The fact is that boric acid is a multibasic acid, i.e., capable of forming multivalent ions during dissociation, which, in turn, can form complexes with base ions (tris). In other words, the scheme of dissociation reactions is significantly different from
the reactions of a monobasic acid and a monobasic base (see.~\eqref{Pt:01.15}).

In~\cite{Michov} it is indicated that in solution boric acid $HB$ and tris $H^+T$, in addition to the acid residue $B^-$ and the main residue $T$, forms complexes are the Tris-boric acid $HTB$,  three-boric acid $H_3B_3$, triboration $H_2B_3^-$, Tris-borate ion $TB^-$.

\subsection{Scheme of reactions}\label{Pt:sec:03.01}

In principle, there are different schemes of chemical reactions in a solution of boric acid ($H_3BO_3$) + tris ($C_4H_{11}NO_3$). Here the variant proposed in~\cite{Michov} (see also \cite[p.\, 64--66]{PolyakovaZhukovShiryaeva}) is used the notations used in~\cite{Michov}.
\begin{equation}\label{Pt:03.01}
BH   \overset{pK_1=9.29}  \rightleftarrows   B^-+H^+,
\end{equation}
\begin{equation}\label{Pt:03.02}
3HB  \overset{pK_2=-1.77} \rightleftarrows   H_3B_3,
\end{equation}
\begin{equation}\label{Pt:03.03}
HB_3 \overset{pK_3=9.02}  \rightleftarrows   H_2B_3^-+H^+,
\end{equation}
\begin{equation}\label{Pt:03.04}
HB+T  \overset{pK_4=-2.53} \rightleftarrows   HTB,
\end{equation}
\begin{equation}\label{Pt:03.05}
HTB    \overset{pK_5=9.50} \rightleftarrows   TB^-+H^+,
\end{equation}
\begin{equation}\label{Pt:03.06}
H^+T   \overset{pK_6=7.98} \rightleftarrows   T+H^+,
\end{equation}
\begin{equation}\label{Pt:03.07}
H_2O   \overset{pK_w=14.0} \rightleftarrows   H^++OH^-,
\end{equation}
where $HB$ is the boric acid, $H^+T$ is the tris (more precisely, a single-charge base ion), $B^-$ is the a single-charge boric acid ion, $T$ is the main residue, $HTB$ is the tris-boric acid,  $H_3B_3$ is the three-boric acid, $H_2B_3^-$ is the triboration,  $TB^-$ is the tris-borate ion, $H^+$ is the hydrogen ion, $OH^-$ is the hydroxyl ion, the numbers above the arrows indicate the dissociation constants of the reactions on a logarithmic scale and are related to $K_i$ by relations
\begin{equation}\label{Pt:03.08}
K_i = 10^{-pK_i}, \quad K^2_w = 10^{-pK_w}.
\end{equation}

The equation \eqref{Pt:03.07} is the equation of water autoprotolysis. It has no direct relation to the dissociation reactions of boric acid and tris. The inclusion of the equation \eqref{Pt:03.07} in the general scheme of reactions indicates only that the contribution of the ions $OH^-$ and $H^+$ will be taken into account in the equation of electroneutrality.

It is appropriate to note that according to Arrhenius (see e.g., \cite[p.\, 63]{Korita}) an acid $HA$ is defined as a substance that cleaves hydrogen ions in solution
\begin{equation}\label{Pt:03.09}
HA \overset{pK_a} \rightleftarrows   A^- +H^+,
\end{equation}
and the base $BOH$ is defined as a substance that cleaves hydroxyl ions in solution
\begin{equation}\label{Pt:03.10}
BOH  \overset{pK'_b} \rightleftarrows   B+ + OH^-,
\end{equation}
where $pK_a$, $pK'_b$ are the equilibrium constants of acid and base (here the stroke symbol does not mean taking into account the ionic strength).

The equilibrium equations for the reactions \eqref{Pt:03.08}, \eqref{Pt:03.09} are written in the form
\begin{equation}\label{Pt:03.11}
\dfrac{[A^-][H^+]}{[HA]} = K_a, \quad \dfrac{[B^+][OH^-]}{[BOH]} = K'_b.
\end{equation}
In addition, from the equilibrium condition for the equation \eqref{Pt:03.07} we have (autoprotolysis of water)
\begin{equation}\label{Pt:03.12}
[OH-][H^+] = K_w^2 \quad \Rightarrow \quad [OH^-]=\dfrac{K^2_w}{[H^+]}.
\end{equation}
This means, in particular, that the second equation \eqref{Pt:03.11} can be written in the form
\begin{equation}\label{Pt:03.13}
\dfrac{[B^+] K^2_w}{[H^+][BOH]} = K'_b  \quad \text{or} \quad \dfrac{[H^+][BOH]}{[B^+]} = \dfrac{K^2_w}{K'_b} = K_b
\end{equation}
This equilibrium corresponds to the reaction
\begin{equation}\label{Pt:03.14}
B^+ \overset{pK_b} \rightleftarrows   BOH + H^+ \quad  \Leftrightarrow \quad (H^+B \rightleftarrows B + H^+),
\end{equation}
which formally corresponds to the definition of an acid is the substance $B^+$ cleaves off the ion $H^+$.

From this point of view, bases can be considered as acids (conjugated acids), by formally introducing the notation $B^+ \equiv H^+B$, $BOH \equiv B$. At present time, this is how it is customary to write chemical equations are the consider bases as acids (conjugated), for more information, see e.g., \cite[p.\, 60--87]{Korita}.

\subsection{Molar concentrations}\label{Pt:sec:03.02}

In the general case, when determining the analytical concentrations of $a_s$ for six ($r=6$) chemical equations \eqref{Pt:03.01}--\eqref{Pt:03.06} we should write the matrix of stoichiometric coefficients $\nu_{ik}$, $i=1,\dots,r$, $k=0,\dots,n$, where $n$ is the number of components of the mixture excluding the hydrogen ion (in this case, $n=8$, i.e., the components are $HB$, $H^+T$, $B^-$, $T$, $HTB$, $H_3B_3$, $H_2B_3^-$, $TB^-$). Then the fundamental solutions of the system are determined
\eqref{Pt:02.12}, i.e., $\lambda_k^{(s)}$, $s=1,\dots,r_0$, $k=0,\dots,n$, where $r_0 = \rank \nu_{ik}$, and finally, using the relations \eqref{Pt:02.13} the analytical concentrations, i.e., the values $a_s$, $s=1,\dots,r_0$ are found.

In the case of using the chemical form of recording (i.e., as in this case, using the symbols of chemical elements) the procedure for finding analytical concentrations is significantly simplified and is possible without tedious writing out the matrix of stoichiometric coefficients $\nu_{ik}$ and without finding the fundamental solutions $\lambda_k^{(s)}$. The following set of rules should be used.

1. The rank of the matrix $r_0$ is equal to the number of initial substances before dissociation. The number of analytical concentrations is $r_0$, in this case, $r_0 = 2$, there are two initial substances are the boric acid and the tris.

2. The analytical concentration is a linear combination of the concentrations of all substances that have the same chemical symbol. For example, the analytical concentration of boric acid (symbol $B$) will include the concentrations of substances containing the symbol $B$.

3. If the chemical symbol contains a lower index, for example, $B_k$, then the concentration of the corresponding component is included in the analytical concentration in the form of $k[B_k]$.

Using the above rules for the molar analytical concentrations of boric acid and tris, we obtain
\begin{equation}\label{Pt:03.15}
 C_T=[T]+[H^+T]+[TB^-]+[HTB],
\end{equation}
\begin{equation}\label{Pt:03.16}
 C_B=[HB]+[B^-]+3[H_3B_3]+3[H_2B_3^-]+[TB^-]+[HTB],
\end{equation}
where $C_B$, $C_T$ are the molar analytical concentrations of boric acid and tris.

The results obtained in par.~\ref{Pt:sec:03.02}, correspond to Step 1 of the algorithm (see~p.~\pageref{Pt:etap:01}).

\subsection{Equilibrium equations of reactions}\label{Pt:sec:03.03}

In the case of reactions \eqref{Pt:03.01}--\eqref{Pt:03.06} the equilibrium equations have the form
\begin{equation}\label{Pt:03.17}
\frac{[B^-][H^+]}{[HB]} = K_1 = 10^{-pK_1}, \quad pK_1 = 9.29,
\end{equation}
\begin{equation}\label{Pt:03.18}
\frac{[H_3B_3]}{[HB]^3} = K_2 = 10^{-pK_2}, \quad pK_2 = -1.77,
\end{equation}
\begin{equation}\label{Pt:03.19}
\frac{[H_2B_3^-][H^+]}{[H_3B_3]} = K_3 = 10^{-pK_3}, \quad pK_3 = 9.02,
\end{equation}
\begin{equation}\label{Pt:03.20}
\frac{[HTB]}{[HB][T]} = K_4 = 10^{-pK_4}, \quad pK_4 = -2.53,
\end{equation}
\begin{equation}\label{Pt:03.21}
\frac{[TB^-][H^+]}{[HTB]} = K_5 = 10^{-pK_5}, \quad pK_5 = 9.50,
\end{equation}
\begin{equation}\label{Pt:03.22}
\frac{[T][H^+]}{[H^+T]} = K_6 = 10^{-pK_6}, \quad pK_6 = 7.98.
\end{equation}

Recall that to solve the problem, i.e., to determine the dependence of the unknowns $[HB]$, $[H^+T]$, $[B^-]$, $[T]$, $[HTB]$, $[H_3B_3]$, $[H_2B_3^-]$, $[TB^-]$ ($n=8$) on $[H^+]$, $C_B$, $C_T$ of the form \eqref{Pt:02.16} we have $8$ of nonlinear equations are the $6$ of equations \eqref{Pt:03.17}--\eqref{Pt:03.22} and the $2$ equations \eqref{Pt:03.15}, \eqref{Pt:03.16}  (of course, this is the case when $I=0$).

In accordance with Step 2 of the algorithm (see~p.~\pageref{Pt:etap:02}) it is necessary to plot the dependence of the concentration of $[HB]$, $[H^+T]$, $[B^-]$, $[T]$, $[HTB]$, $[H_3B_3]$, $[H_2B_3^-]$, $[TB^-]$ on the values of $[H^+]$, $C_B$, $C_T$. However, as the analysis of the equations shows, it is not possible to obtain an explicit form of such a dependence.
According to the Remark $2$ (see~p.~\pageref{Pt:note:02}), instead of $[H^+]$, $C_B$, $C_T$ some other set consisting of $\xi_0$ and $r_0$ ‘parameters’ must be chosen.

In this considered case, using \eqref{Pt:03.17}--\eqref{Pt:03.22} it is convenient to write the dependencies of the concentrations $[H_2B_3^-]$,
$[TB^-]$, $[B^-]$, $[H^+T]$, $[H_3B_3]$, $[HTB]$ from concentrations $[HB]$, $[T]$ (and $[H^+]$)
\begin{equation}\label{Pt:03.23}
[H_2B_3^-] = \frac{[HB]^3}{[H^+]} 10^{- pK_3 - pK_2} = \frac{[HB]^3}{[H^+]} K_3 K_2,
\end{equation}
\begin{equation}\label{Pt:03.24}
[TB^-] = \frac{[HB][T]}{[H^+]} 10^{- pK_5 - pK_4} = \frac{[HB][T]}{[H^+]} K_5K_4,
\end{equation}
\begin{equation}\label{Pt:03.25}
[B^-] = \frac{[HB]}{[H^+]} 10^{- pK_1} = \frac{[HB]}{[H^+]} K_1,
\end{equation}
\begin{equation}\label{Pt:03.26}
[H^+T] = [T][H^+] 10^{+ pK_6} = \dfrac{[T][H^+]}{K_6},
\end{equation}
\begin{equation}\label{Pt:03.27}
[H_3B_3] = [HB]^3 10^{- pK_2} = [HB]^3 K_2,
\end{equation}
\begin{equation}\label{Pt:03.28}
[HTB] = [HB][T] 10^{- pK_4} = [HB][T] K_4.
\end{equation}
Of course, in this case the analytical concentrations $C_B$, $C_T$ must also be expressed in terms of $[HB]$, $[T]$ (and $[H^+]$) using the relations \eqref{Pt:03.15}, \eqref{Pt:03.16}.

The results obtained in~par.~\ref{Pt:sec:03.03}, correspond to Step 2 of the algorithm (see~par.~\pageref{Pt:etap:02}). The number of nonlinear equations that will need to be solved has been reduced to $r_0 + 1 = n - r + 1$ (in this case $r_0 + 1 = 3$, for $I=0$ or $r_0 + 2 = 4$, for $I \ne 0$). Recall that initially for $I \ne 0$ there were $n + 2 =10$ equations.

\subsection{Determination of $pH$ of the solution}\label{Pt:sec:03.04}

To determine the concentrations of $[HB]$, $[T]$, $[H^+]$ to the equations \eqref{Pt:03.15}, \eqref{Pt:03.16}
\begin{equation}\label{Pt:03.29}
 C_T=[T]+[H^+T]+[TB^-]+[HTB],
\end{equation}
\begin{equation}\label{Pt:03.30}
 C_B=[HB]+[B^-]+3[H_3B_3]+3[H_2B_3^-]+[TB^-]+[HTB],
\end{equation}
adding the electroneutrality equation
\begin{equation}\label{Pt:03.31}
   [H^+T] - [H_2B_3^-] - [TB^-] - [B^-] + [H^+] - \frac{K_w^2}{[H^+]}=0.
\end{equation}
From the equations \eqref{Pt:03.29}--\eqref{Pt:03.31} using \eqref{Pt:03.23}--\eqref{Pt:03.28} we exclude the quantities $[H_2B_3^-]$,
$[TB^-]$, $[B^-]$, $[H^+T]$, $[H_3B_3]$, $[HTB]$.

Then
\begin{equation}\label{Pt:03.32}
 C_T=[T]+\dfrac{[T][H^+]}{K_6}+ K_4K_5\frac{[HB][T]}{[H^+]}+ K_4[HB][T],
\end{equation}
\begin{equation}\label{Pt:03.33}
C_B=[HB]+ K_1\frac{[HB]}{[H^+]}+3K_2[HB]^3 + 3 K_2 K_3\frac{[HB]^3}{[H^+]} + {}
\end{equation}
\begin{equation*}
{} +  K_4K_5\frac{[HB][T]}{[H^+]}+ K_4[HB][T],
\end{equation*}
\begin{equation}\label{Pt:03.34}
   \dfrac{[T][H^+]}{K_6} -  K_2 K_3\frac{[HB]^3}{[H^+]} - K_4K_5 \frac{[HB][T]}{[H^+]} -  K_1\frac{[HB]}{[H^+]} + [H^+] - \frac{K_w^2}{[H^+]}=0.
\end{equation}
Note that the full version of the electroneutrality equation is written here. In principle, in some cases, approximate variants can be used, neglecting the concentration of $H^+$ and/or $OH^-$ ions.

In the physical sense, the concentrations of $[HB]$, $[T]$, $[H^+]$ must be positive (it is easy to check that there are no zero solutions). In view of this, we should make substitutions for similar  \eqref{Pt:02.11}, e.g.,
\begin{equation}\label{Pt:03.35}
   [HB] = 10^X, \quad [T] = 10^Y, \quad [H^+] = 10^{-pH}
\end{equation}
and solve the system  \eqref{Pt:03.32}--\eqref{Pt:03.34}  with respect to $X$, $Y$, $pH$. Note that, in principle, any base can be chosen (not necessarily $10$).
The base $10$ is chosen because of tradition, since the determination of $pH$ and the dissociation constants $pK_i$ is made in the chemistry literature using the decimal logarithm $\lg$.

\subsubsection{Correction scheme}\label{Pt:sec:03.04.01}

Recall that the system \eqref{Pt:03.32}--\eqref{Pt:03.34} is written for the case when the ionic strength $I=0$. In order to take into account the influence of the ionic strength, the dissociation constants $K_i$ should be replaced by $K'_i$ using the relations \eqref{Pt:02.04}--\eqref{Pt:02.06}. In the case under consideration, all the ions are singly-charged. In particular, this means that all the activity coefficients for charged ions will be the same
\begin{equation}\label{Pt:03.36}
  \lg \gamma = \lg \gamma_k = - A z_k^2 I^{1/2}.
\end{equation}
The formulas \eqref{Pt:02.04} for the constants $K'_i$ take the form
\begin{equation}\label{Pt:03.37}
  K'_1 = K_1 \gamma^{-2}, \quad
  K'_2 = K_2 \gamma^{3}, \quad
  K'_3 = K_3 \gamma^{-2},
\end{equation}
\begin{equation*}
  K'_4 = K_4, \quad
  K'_5 = K_5 \gamma^{-2}, \quad
  K'_6 = K_6.
\end{equation*}
Note that it is not necessary to use the formula \eqref{Pt:02.04} at all. It is enough to write down the equations \eqref{Pt:03.17}--\eqref{Pt:03.22} for the activities (see \eqref{Pt:01.30}, \eqref{Pt:01.31}), e.g.,
\begin{equation}\label{Pt:03.38}
\frac{[B^-]\gamma_{B^-}[H^+]\gamma_{H^-}}{[HB] \gamma_{HB}} = K_1 \quad \Rightarrow \quad K'_1 = \dfrac{\gamma_{HB}}{\gamma_{B^-}\gamma_{H^-}}K_1,
\end{equation}
\begin{equation*}
\gamma_{B^-} = \gamma, \quad \gamma_{H^-} = \gamma, \quad \gamma_{HB}= 1.
\end{equation*}

As shown by the computational experiment, the direct substitution of  \eqref{Pt:03.37} into the equations \eqref{Pt:03.32}--\eqref{Pt:03.34} significantly complicates the nonlinear system of equations. The numerical solution algorithm (Newton's method) has a small convergence rate and requires a very good accurate initial approximation. In contrast, the algorithm for solving the system \eqref{Pt:03.32}--\eqref{Pt:03.34} with constant coefficients $K_i$ (the same Newton's method) has fast convergence and converges from almost any initial approximation.

In view of the above, an iterative correction algorithm was used.

1. The solution of the system  \eqref{Pt:03.32}--\eqref{Pt:03.34} with constant coefficients $K_i$ is calculated.

2. For the obtained solution, the ionic strength $I$ is calculated and the constants are corrected using \eqref{Pt:03.38}.

3. The convergence condition is checked (by the value of $pH$) and if it is not met, it returns to par. 1.

\textbf{Remark 3.} In particular, in~\cite[p.\, 93]{Korita}  states that instead of calculating $pH$ using the formula $pH=-\lg[H^+]$  the activity of the hydrogen ion \label{Pt:note:03} should be used instead of the concentration
\begin{equation}\label{Pt:03.39}
pH_a = -\lg(\gamma_{H^+}[H^+]).
\end{equation}

\label{Pt:note:04}
\textbf{Remark 4.} Instead of the relations \eqref{Pt:02.05}, \eqref{Pt:02.06} more precise formulas can be used (for more information, see e.g., in~\cite[pp.\, 38--48]{Korita}). This does not affect the convergence of the solution process, at least for the chemical equations under consideration.

\newpage
\setcounter{equation}{0}
\setcounter{figure}{0}

\section{Calculation results}\label{Sh:sec:04}

Here are the results of calculations for chemical reactions in the case of a Tris-borate buffer mixture \eqref{Pt:03.01}--\eqref{Pt:03.06}, using the results of par.~\ref{Pt:sec:02}, \ref{Pt:sec:03}. The dissociation constants specified in the formulas were chosen as parameters \eqref{Pt:03.17}--\eqref{Pt:03.22} (see~\cite{Michov})
\begin{equation}\label{Pt:04.01}
pK_1 = 9.29, \quad
pK_2 = -1.77, \quad
pK_3 = 9.02,
\end{equation}
\begin{equation*}
pK_4 = -2.53, \quad
pK_5 = 9.50, \quad
pK_6 = 7.98, \quad
pK_w = 14.00,
\end{equation*}
which are related to $K_i$ by the relations
\begin{equation}\label{Pt:04.02}
K_i = 10^{-pK_i}.
\end{equation}

The molar analytical concentrations of boric acid $C_B$ and tris $C_T$ were chosen from the intervals determined by the ratios
\begin{equation}\label{Pt:04.03}
0. 1 \leq C_B \leq 0.3, \quad 0.1 \leq C_T \leq 0.3.
\end{equation}

The system of nonlinear equations \eqref{Pt:03.32}--\eqref{Pt:03.34} with respect to variables $X$, $Y$, $pH$ (see \eqref{Pt:03.35}) for constants $K_i$ was solved by Newton's method (see, e.g., \cite{Bahvalov,Kalitkin}) with subsequent correction (in the case of a non-zero ionic strength $I \ne 0$), described in~par.~\ref{Pt:sec:03.04.01}. Both standard algorithms  (\textrm{fsolve Maple}), and proprietary original programs were used. The latter, i.e., native programs, were used to solve the system \eqref{Pt:03.32}--\eqref{Pt:03.34}, in which ‘constants’ $K'=K'_i(I)$. As already mentioned, the direct solution of the system \eqref{Pt:03.32}--\eqref{Pt:03.34} with variable coefficients proved ineffective, which led to the use of the correction algorithm.

The solution of the system \eqref{Pt:03.32}--\eqref{Pt:03.34} was calculated with a relative error of $\varepsilon = 10^{-6}$. The same error was set when using the correction mechanism. The calculation results for the concentrations $C_T = 0.2$, $0. 1 \leq C_B \leq 0.3$
are given in Table \ref{Pt:table:02}. The table contains various values of $pH$ are the calculated by the classical formula $pH=-\lg[H^+]$,
by the formula \eqref{Pt:03.39} with taking into account the activity of $pH_a$ and $pH_{I=0}$
without taking into account the influence of the ionic strength. A large number of signs are retained in the table to emphasize the difference between the values. In practice, for the values of $pH$, it is sufficient to keep two or three decimal places, since the accuracy of measuring instruments currently does not exceed $0.01$ units of $pH$. In addition, the table shows the ionic strength $I$ and activity coefficients $\gamma$. The value of $\lg\gamma$ is, in fact, the difference between the values of $pH$ and $pH_a$.

\begin{table}[H]
  \centering

\caption{Dependence of $pH$ on the concentration of boric acid $C_B$ at the concentration of tris $C_T =0.2$}\label{Pt:table:02}

\begin{tabular}{|c|c|c|c|c|c|c|c|}
  \hline
$C_B$& $pH$& $pH_a$& $pH_{I=0}$& $I\cdot 10^{2}$& $\gamma$ & $(\lg\gamma)\cdot 10^{2}$ \\
\hline
0.10& 8.6778 & 8.7448& 8.7444 & 1.731 & 0.8570 & -6.700 \\
0.12& 8.5785 & 8.6453& 8.5783 & 1.719 & 0.8574 & -6.677 \\
0.14& 8.4782 & 8.5435& 8.4766 & 1.641 & 0.8605 & -6.524 \\
0.16& 8.3735 & 8.4361& 8.3706 & 1.510 & 0.8657 & -6.259 \\
0.18& 8.2651 & 8.3244& 8.2614 & 1.353 & 0.8724 & -5.924 \\
0.20& 8.1620 & 8.2180& 8.1583 & 1.210 & 0.8789 & -5.603 \\
0.22& 8.0748 & 8.1284& 8.0719 & 1.109 & 0.8838 & -5.364 \\
0.24& 8.0058 & 8.0579& 8.0039 & 1.047 & 0.8869 & -5.211 \\
0.26& 7.9512 & 8.0024& 7.9500 & 1.010 & 0.8887 & -5.119 \\
0.28& 7.9067 & 7.9574& 7.9060 & 0.989 & 0.8899 & -5.065 \\
0.30& 7.8693 & 7.9196& 7.8689 & 0.977 & 0.8905 & -5.034 \\
\hline
\end{tabular}

\end{table}

From the table~\ref{Pt:table:02} it is clearly seen that the values of  $pH$ and $pH_{I=0}$ differ only by the third decimal point (except for $C_B=0.1$, $C_T=0.2$). In other words, the influence of the ionic strength is insignificant. In contrast, the values of $pH$ and $pH_a$ differ by the second decimal place, which is within the accuracy of the $pH$ measurements. The same pattern is observed for other values of concentrations varying in the intervals indicated by the inequalities \eqref{Pt:04.03}. Note that the relationship between the ionic strength and the activity coefficients was chosen based on the relations \eqref{Pt:02.05}, \eqref{Pt:02.06}
for $A = 0.5093 ~\text{mol}^{-1/2}\text{L}^{1/2}$ (see \eqref{Pt:01.28}). As already mentioned (see Remark 4 on p.~\pageref{Pt:note:04}), instead of
\eqref{Pt:02.05}, \eqref{Pt:02.06}  other ratios may also be used.

\begin{figure}[H]
 \centering\includegraphics[scale=0.95]{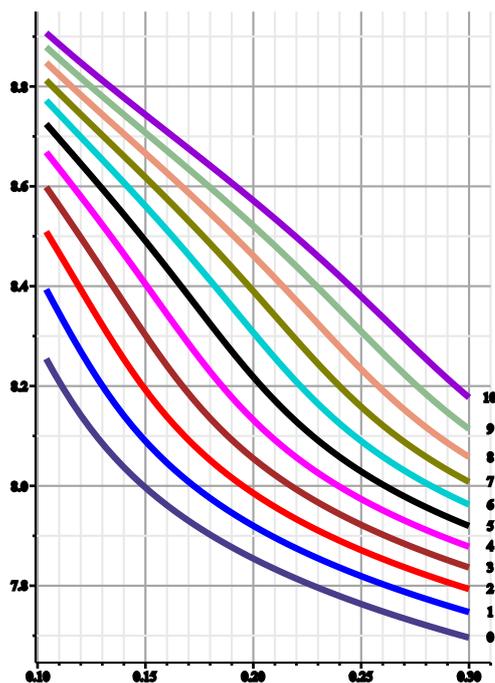}\\
 \caption{The dependencies of $pH_a$ on concentration $C_B$ at fixed values of concentration $C_T$. The curves marked with numbers $n=0,\dots,10$, correspond to concentrations $C_T = 0.1 + 0.02\,n$}
\label{Pt:fig:4.01}
\end{figure}

Figure~\ref{Pt:fig:4.01} shows the results of calculations of the value $pH_a$ for the intervals of change $C_B$, $C_T$, indicated by the inequalities \eqref{Pt:04.03}

Analysis of the results shows that the ‘Tris-borate’ mixture can be used to obtain $pH_a$ solution in the interval
\begin{equation}\label{Pt:04.04}
7.7 < pH_a < 8.9, \quad (0. 1 \leqq C_B \leqq 0.3, \quad 0.1 \leqq C_T \leqq 0.3).
\end{equation}

When $C_B$, $C_T$ values are fixed, a stable value of $pH_a$ can be created, e.g., for the zonal electrophoresis process. In the case when, e.g., $C_T$ is fixed, by changing the value of  $C_B$ along the electrophoretic chamber can create a $pH_a$ distribution (the so-called $pH$-gradient) and use a Tris-borate buffer mixture to conduct the isoelectric focusing process (see par.~\ref{Pt:sec:01}).

Quite often, the electrophoretic chamber in which isoelectric focusing is performed is a long cylindrical region or a long flat gel plate, which allows us to consider the region as one-dimensional (see e.g., \cite{Righetti1990,Righetti1986,Righetti1986-1}). Setting
concentration distribution  $C_B=C_B(x)$ and knowing the dependence $pH_a=pH_a(C_B)$,  it is easy to get the dependence $pH_a(x)=pH_a(C_B(x))$. Of particular interest are the linear $pH_a(x)$,  the so-called linear $pH$-gradients. The fact is that, in particular, with a linear distribution of $pH_a(x)$ it is convenient to interpret the results of the mixture separation process.

The computational experiment showed that the dependence $pH_a=pH_a(C_B)$ is practically linear at $C_T \geqq 0.3$ (see, in particular, $10$ line in Figure~\ref{Pt:fig:4.01}). Using the least-squares method (e.g., \textrm{Fit Maple}), it is easy to get a polynomial dependence (in this case linear)
\begin{equation}\label{Pt:04.05}
pH_a \approx  9.299 - 3.694\,C_B, \quad C_T=0.3.
\end{equation}
The fitting is held at $100$ points. The standard deviation is $\sigma = 0.0084$,  and the maximum (point) deviation from the exact formula is  $max_0 = 0.015$.

Using a quadratic approximation improves the result
\begin{equation}\label{Pt:04.06}
pH_a \approx 9.206 - 2.672\, (C_B) - 2.555\,(C_B)^2, \quad C_T=0.3,
\end{equation}
giving the deviation $max_0 = 0.0097$.

To perform the calculations, the created program (\textrm{Maple}) was used, the listing of which with detailed comments is given in the Appendix.

\newpage
\newpage
\addcontentsline{toc}{section}{Conclusion}
\section*{Conclusion}

Stationary equilibria were found for the equilibrium chemical reactions of the complex ‘Tris-borate’ mixture and the influence of the ionic strength of the solution on $pH$ of the mixture was investigated. The last stage of the work, which made it possible to obtain numerical results, required considerable preliminary work. The system of nonlinear equations consisting of six equilibrium equations and the equation of electroneutrality for the eight components of mixtures formed as a result of dissociation was converted to three nonlinear equations, which significantly simplified the numerical solution. It turned out that to take into account for the effect of the ionic strength of the solution on the $pH$ solution, it is convenient to solve the problem in two stages. At the first stage the problem should be solved with constant dissociation constants (i.e., without taking into account the influence of the ionic strength on the dissociation constants), and at the second stage the solution should be corrected taking into account the ionic strength of the solution.

From the above results, it is clear that a simple solution to the problem of calculating chemical equilibria is possible only for relatively elementary equations of chemical kinetics. In the case when the occurrence of various chemical complexes in the reactions is taken into account, the problem becomes significantly more complicated. Then, without preliminary analysis and simplifications, it is quite difficult to construct direct the solution of the original system of nonlinear equations describing the equilibrium.

The proposed correction algorithm for taking into account the ionic strength turned out to be quite effective, which allowed us to conduct an extensive computational experiment, only a small part of which is presented in the work. The analysis showed that the ‘Tris-borate’ mixture in a certain range of $pH$ values can be used to create a linear distribution of $pH$ solution and conduct the process of isoelectrophocusing.

\newpage

\addcontentsline{toc}{section}{References}
\renewcommand{\refname}{\centering

\end{document}